\documentclass{amsart}
\usepackage{amsmath}
\usepackage{amsfonts}
\usepackage{mathrsfs}
\usepackage{bbm}

\usepackage[normalem]{ulem}

\newcommand{\bpf}[1][Proof]{{\noindent {\sc #1: }}}
\newcommand{\epf}{{{\hfill $\Box$ \smallskip}}}

\newcommand{\R}{\mathbb{R}}
\newcommand{\N}{\mathbb{N}}
\newcommand{\T}{\mathbb{T}}
\newcommand{\Z}{\mathbb{Z}}
\newcommand{\tbf}{\mathbf{t}}
\newcommand{\sbf}{\mathbf{s}}
\newcommand{\ibf}{\mathbf{i}}

\newcommand{\rbf}{\mathbf{r}}
\newcommand{\Pp}{\mathsf{P}}
\newcommand{\Qq}{\mathsf{Q}}
\newcommand{\E}{\mathsf{E}}
\newcommand{\Ic}{\mathcal{I}}
\newcommand{\Bc}{\mathcal{B}}
\newcommand{\wb}{\overline}
\newcommand{\Pc}{\mathcal{P}}
\newcommand{\supp}{\mathop{\rm supp}}

\newtheorem{theorem}{Theorem}
\newtheorem{lemma}{Lemma}
\newtheorem{remark}{Remark}

\begin{document}

\title{Invariant Densities for Dynamical Systems with Random Switching}
\author{Yuri Bakhtin \and Tobias Hurth}

\begin{abstract}
We consider a non-autonomous ordinary differential equation on a smooth
manifold, with right-hand side that randomly
switches between the elements of a finite family of smooth vector fields. For the resulting random dynamical system,
we show that H\"ormander type hypoellipticity conditions are sufficient for uniqueness
and absolute continuity of an invariant measure. 
\end{abstract}

\maketitle



\section{Introduction}

In this paper we study the ergodic theory of systems with random switchings. Such a system
can be described in terms of a finite family of vector fields. We assume that at any given time the
evolution is driven by one of these vector fields, and at random times the driving vector field changes to another
one from the same family. Systems of this nature arise naturally in applications and we refer to the recent
monograph~\cite{Yin} for motivation and extensive bibliography. 

 Many long-term asymptotic properties of dynamical systems or random dynamical systems can be
described in terms of invariant distributions.
The existence of invariant measures often can be derived using the Lyapunov function technique that helps to
establish recurrence properties or tightness, see, e.g.,\cite[Sections~3.3--3.4]{Yin}.

The uniqueness and absolute
continuity of invariant distributions
are often related to each other and more subtle, especially in the case that we consider in this paper where no
diffusion is involved, and the only source of randomness is the random sequence of driving vector fields. Although
some claims have been scattered through the literature, no general result is known, see, e.g.,~\cite[Section
8.5.2]{Yin}.

The goal of this paper is to close this gap and obtain new general conditions that guarantee uniqueness and absolute
continuity of invariant measures for systems with random switchings. The two conditions that we suggest are
formulated in terms of Lie algebras associated to the driving vector fields. They are close analogues of the
classical H\"ormander condition guaranteeing absolute continuity of transition densities of hypoelliptic diffusions.
In the diffusion context, this result is usually derived from the variational analysis of diffusion paths known as
Malliavin calculus, see, e.g.,~\cite[Chapter VIII]{Bass:MR1483890},\cite{Bell:MR2250060},\cite{Nualart:MR2200233}.

In fact, the central part of this paper is the analysis of transition probabilities of switched systems. Under the
first of our conditions, we prove that all transition probabilities for the system have nontrivial absolutely
continuous components. The second condition is more general, and it allows to prove the existence of absolutely
continuous components not for the transition probabilities themselves, but for their time averages. The extraction of
these absolutely continuous components is largely based on classical control theory results that can be found  in
Chapter~$3$ of~\cite{Jurdjevic:MR1425878}. These control theory results rely on earlier work by
Chow~\cite{Chow:MR0001880}, Sussmann and Jurdjevic~\cite{Sussmann-Jurdjevic:MR0338882}, and
Krener~\cite{Krener:MR0383206}. 
Our conditions and the
structure of our proofs match those of~\cite{Jurdjevic:MR1425878}, where the nondegeneracy of certain
maps is exploited to establish the accessibility property. We use the same nondegeneracy to prove absolute continuity,
and one can
interpret our result as filling the control theory with probabilistic content. In fact, the idea to use the
geometric control theory approach to establish regularity of Markov transition kernels along with ergodic
properties is not new, see, e.g., \cite{Agrachev-Kuksin-Sarychev-Shirikyan:MR2329509} where controllability
of the 2D Navier--Stokes system was used to prove the absolute continuity of finite-dimensional projections of
transition kernels.

The article is organized as follows: In Section~\ref{sec:definitions} we introduce the setting, necessary notation and
notions from differential geometry and geometric control theory. We also state the main result on uniqueness and
absolute continuity of invariant measures, and two central auxiliary results on regularity of transition probabilities
each based on one of the H\"ormander type assumptions. We prove these regularity results in
Sections~\ref{sec:AC-component-1-proof} and~\ref{sec:AC-component-2-proof}. In
Section~\ref{sec:AC-of-ergodic-measures} we prove that any ergodic measure has to be absolutely continuous if
its support contains a point where hypoellipticity holds. Section~\ref{sec:uniqueness-proof} contains the
proof of the main result: if the hypoellipticity holds at a point that can be approached from any initial point using
the given vector fields as admissible controls, then there exists at most one invariant distribution, and this
distribution has to be absolutely continuous. In Section~\ref{sec:examples}, we apply the main result to a switching system on the $n$-dimenisonal torus and a switching system involving two Lorenz vector fields.

{\bf Acknowledgments:} The idea to study invariant densities for systems with switchings emerged after a discussion
of ``blinking systems'' with Leonid Bunimovich and Igor Belykh, and we would like to thank them. We are also thankful to
L.Bunimovich for his comments on the Lorenz system.
We are grateful
to Martin Hairer and especially Jonathan Mattingly for stimulating discussions of other possible approaches to the main
results of this paper.  We thank the referees for their comments that helped us to improve the paper. The partial
support from NSF through a CAREER Award DMS-0742424 is gratefully acknowledged by
YB.
\section{Definitions, Notation, and Main Results} \label{sec:definitions}

\subsection{The dynamics} \label{sec:dynamics}

We consider a finite collection $D$ of smooth and forward complete vector fields 
on an $n$-dimensional $C^\infty$-manifold $M$. 

We denote these vector fields by $u_i,i\in S=
\{1,\ldots,k\}$.
Each vector field $u$ in $D$ induces an ordinary differential equation of the
form 
\begin{equation*}
\dot{x}(t) = u(x(t)).
\end{equation*}
This differential equation is uniquely solvable if equipped with
an initial condition
\begin{equation*}
x(0) = \xi\in M,
\end{equation*}
 and forward completeness means that the solution trajectories are well-defined
for all times $t>0$. 

We can define a stochastic process
$X=(X_t)_{t \geq 0}$ on $M$ in the following way: Given an initial state $i\in S$ and an initial
value $\xi \in M$, $X_t$ follows the trajectory generated by the vector field $u_i$  and initial
condition $\xi$ for an 
exponentially distributed random time with parameter $\lambda_i>0$. Then a new state is selected at random from
$S\setminus\{i\}$, and,
for another exponentially distributed random time, $X_t$ follows the new vector field corresponding to that state. 
Iterating this construction  we obtain a piecewise smooth
trajectory $(X_t)_{t\ge 0}$ defined for all positive times and driven by one of the vector fields from $D$
between any two switchings.
We assume that (i) all the inter-switching times are exponentially distributed and
independent conditioned on the sequence of driving vector
fields, (ii) the parameter $\lambda_j$ of the exponential time between any two switches depends only on the current
state $j$, 
 and (iii) the probabilities of
switchings between any two states are positive.

We choose to work with exponential
waiting times to ensure the Markov property, although our results can be extended to
non-Markovian settings resulting from more general waiting time distributions.

It is convenient to keep track of the driving vector fields at all times. We
define $A_t\in S$ as the index of the driving vector field at time $t$, also
referred to as the regime or state at time $t$. It is a Markov process with continuous time and finitely
many states. Its trajectories are right-continuous and piecewise constant.

Although $X$ alone is not a Markov process, the joint process $(X,A)$ is Markov. 
We denote elements of the
associated Markov family, i.e., the distribution on paths emitted at $(\xi,i)\in M \times S$  and generated by
the iterative random procedure above, by $\Pp_{\xi,i}$, and the corresponding transition probability measures by
$\Pp_{\xi,i}^{t}$,  $t \geq 0$. The transition probability measures are
defined on the product $\sigma$-algebra $\Bc(M) \otimes \Pc(S)$, where $\Bc(M)$ is the
Borel $\sigma$-algebra on $M$ and $\Pc(S)$ is the power set of $S$. We write $\E_{\xi,i}$ for expectation with
respect to $\Pp_{\xi,i}$.

Let us recall that if the initial distribution of the Markov process $(X,A)$ is $\mu$, then  the distribution of the
process at time $t$ is given by the measure $\mu \Pp^t$ on $M\times S$ defined by
\begin{equation}
 \mu \Pp^t(E\times\{j\}) = \sum_{i=1}^{k} \int_{M} \Pp_{\xi,i}^t(E\times\{j\})\, \mu(d\xi\times\{i\}).
\label{eq:convolution} 
\end{equation}

A probability measure $\mu$ on $M \times S$ is called
invariant for $(\Pp^t)$ if  
$\mu=\mu \Pp^t$ for all $t \geq 0$.

The main goal of this paper is to give conditions on $D$ that would guarantee absolute
continuity and uniqueness of an invariant measure of the Markov semigroup $(\Pp^t)=(\Pp^t)_{t\ge 0}$. The fairly general
conditions that we suggest are formulated in geometric terms, and we proceed to introduce the
necessary
definitions and notation.

\subsection{Auxiliary definitions and notation}\label{sec:diff_geometry}

Let $V(M)$ denote the set of real smooth vector fields on the manifold $M$, and let $C^{\infty}(M)$ denote
the set of real-valued smooth functions on $M$. As explained above, we assume that $D$ is
contained in $V(M)$. Any element of $V(M)$ corresponds uniquely to a derivation on $C^{\infty}(M)$, that is to
a linear operator $\delta$ on $C^{\infty}(M)$ satisfying the Leibniz rule
\begin{equation*}
\delta(f \cdot g) = \delta(f) \cdot g + f \cdot \delta(g).
\end{equation*}    
The Lie bracket of two vector fields $u$ and $v$ in $V(M)$ is defined as the vector field  
\begin{equation*}
[u,v](f):=u(v(f)) - v(u(f))
\end{equation*} 
for test functions $f$ in $C^{\infty}(M)$. The set $V(M)$ equipped with the bilinear operator $[.,.]$ becomes
a Lie algebra over the reals. A subset of $V(M)$ is called involutive if it is closed under taking the Lie bracket. An
involutive subspace of $V(M)$ is called a subalgebra of~$V(M)$.

The smallest subalgebra of $V(M)$ that contains $D$ is denoted $\Ic(D)$. The
derived algebra $\Ic'(D)$ is the smallest algebra containing 
Lie brackets of vector fields in $\Ic(D)$. We have $\Ic'(D)\subset \Ic(D)$, but
$\Ic'(D)$ might not contain any elements of $D$ and may therefore be strictly contained in $\Ic(D)$.
Further, we define $\Ic_{0}(D)$ as the set of vector fields of the form
\begin{equation*}
v + \sum_{i=1}^{k}{\lambda_{i}u_{i}},
\end{equation*}
where $v \in \Ic'(D)$, $u_1,\ldots,u_k \in D$ and $\sum_{i=1}^{k}{\lambda_{i}}=0$. Finally, we set
\begin{equation*}
\Ic(D)(\xi):=\{u(\xi):\ u \in \Ic(D)\}
\end{equation*}
and
\begin{equation*}
\Ic_{0}(D)(\xi):=\{u(\xi):\ u \in \Ic_{0}(D)\}
\end{equation*}
for any $\xi \in M$. The sets $\Ic(D)(\xi)$ and $\Ic_{0}(D)(\xi)$ are finite-dimensional
vector spaces. 

\medskip

Our main results will be based on the following assumptions that can naturally be called hypoellipticity conditions
in analogy with H\"ormander's theory.
We say that a point $\xi\in M$ satisfies Condition~A if $\dim \Ic_{0}(D)(\xi)=n$. 
We say that a point $\xi\in M$
satisfies Condition~B if  $\dim \Ic(D)(\xi)=n$. 

The set of points satisfying Condition~A is open and so is the set of
points satisfying Condition~B.
\medskip

For our absolute continuity results we will need a reference measure on $M$ that will play the role of Lebesgue
measure. As a smooth manifold, $M$ can be endowed with a Riemannian metric. 
The metric tensor can be used to define measures on coordinate patches of $M$. 
One can use then a partition of 
unity in a standard way (see, e.g., \cite[Section 7]{Taylor:MR2245472})
to construct a Borel measure on $M$ whose
pushforward to $\mathbb{R}^{n}$ under any chart map is equivalent to Lebesgue measure. We call the measure on $M$
obtained through this construction Lebesgue measure, denote it by $\lambda^{M}$, and use
it as the main reference measure, often omitting ``with respect to Lebesgue measure'' when writing about absolute
continuity. 
The product of the Lebesgue measure on $M$ and counting measure on $S$ will be called
the Lebesgue measure on
$M\times S$. We denote the Lebesgue
measure on
$\mathbb{R}^{m}$ by~$\lambda^{m}$. 

\medskip

It remains to introduce the flows generated by vector fields in $D$ and the concept of reachability.

For $i\in S$, we denote the flow function of the vector field $u_i$ by $\Phi_i$. Due to forward completeness of $u_i$, the
flow function is uniquely defined for all $t>0$ and $\eta\in M$ by
\begin{align*}
 \frac{d}{dt}\Phi_i(t,\eta)&=u_i(\Phi_i(t,\eta)),\\
 \Phi_i(0,\eta)&=\eta.
\end{align*}

For $m\in\N$, we will consider vectors 
${\bf t}=(t_1,\ldots,t_m)$ of waiting times between subsequent switches and vectors ${\bf i}=(i_1,\ldots,i_m)$ of
driving states during these waiting intervals. We will restrict ourselves to positive waiting times, but it can also be useful 
(see \cite{Sussmann-Jurdjevic:MR0338882} and \cite{Jurdjevic:MR1425878}) to admit flows backwards in
time. 

We write $\R_{+}$ to denote the positive real line $(0;\infty)$. 

For  ${\bf t}=(t_1,\ldots,t_m)\in \R_{+}^{m}$ and ${\bf i}=(i_1,\ldots,i_m)\in S^{m}$, we define 
\begin{equation*}
\Phi_{{\bf i}}({\bf t},\xi):=\Phi_{i_m}(t_m,\Phi_{i_{m-1}}(t_{m-1},\ldots\Phi_{i_1}(t_1,\xi))\ldots)
\end{equation*}
as the cumulative flow along the trajectories of $u_{i_1},\ldots,u_{i_m}$ with starting point $\xi\in M$.

The transition probabilities $\Pp_{\xi,i}^{t}$ can be expressed in terms of cumulative flows. We do not specify
these straightforward relations in order to avoid heavy notation.

A point $\eta\in M$ is called $D$-reachable from a point $\xi\in M$ if there exist a time vector ${\bf
t}$ with
positive components and a vector~${\bf i}$ of driving states such that
\begin{equation*}
\eta = \Phi_{{\bf i}}({\bf t},\xi).
\end{equation*}
If the components of ${\bf t}$ sum up to $t$, we say that $\eta$ is $D$-reachable from $\xi$ at time~$t$.  

For $\xi \in M$ and $t>0$, let $L_{t}(\xi)$ denote the set of $D$-reachable points from $\xi$ at time~$t$, and
let $L(\xi)=\bigcup_{t>0} L_{t}(\xi)$ denote the set of $D$-reachable points from $\xi$.  The points in
the closure $\overline{L(\xi)}$ can be called $D$-approachable from $\xi$.
Let $L=\bigcap_{\xi \in M} \overline{L(\xi)}$ denote the set of points that are $D$-approachable from all
other
points.

\subsection{Main results} \label{sec:absolute-continuity-invariant}

The following is the main theorem of this paper.

\begin{theorem} \label{thm:uniqueness}
Suppose Hypoellipticity Condition~B is satisfied at some $\xi \in L$.
If $(\Pp^t)$ has an invariant measure, then it is unique and absolutely continuous with respect
to the Lebesgue measure on $M\times S$.
\end{theorem}          

\begin{remark}\rm
Of course, Theorem~\ref{thm:uniqueness} remains true if we replace $L$ by any of its subsets. For example, if
one of the vector fields in $D$ has a minimal global
attractor, then it is sufficient to check hypoellipticity for some point of the attractor.
\end{remark}

Uniqueness of invariant distributions is tightly connected to the regularity of the Markov semigroup. Various
aspects of regularity in connection with ergodicity have been studied in the literature: the existence of minorizing
kernels, the strong Feller property, etc. 
The main task in the proof of Theorem~\ref{thm:uniqueness} is to establish regularity for
transition probabilities under Hypoellipticity Condition~B. However, we begin with a much stronger regularity property that 
can be established under the stronger Hypoellipticity Condition~A.

\begin{theorem}\label{thm:AC-component-1}
If Condition~A is satisfied at a point $\xi\in M$, then for any $i \in
S$ and any $t > 0$, the transition kernel $\Pp_{\xi,i}^{t}$ has a nonzero absolutely
continuous component with respect to Lebesgue measure on $M \times S$. 
\end{theorem}

Under the weaker Condition~B it may happen that none of the transition 
probability measures $\Pp_{\xi,i}^t$, $t > 0$, has a nonzero absolutely continuous component. For example, let $M$~be
the
$n$-dimensional torus $\T^n=\R^n/\Z^n$, and let 
$D=\{u_1,\ldots,u_n\}$ be the standard basis in $\mathbb{R}^{n}$. Fix an arbitrary time $t > 0$. The set of
points $D$-reachable from the origin at time $t$ is the image of 
\begin{equation*}
\biggl\{(s_1,\ldots,s_n) \in [0;\infty)^{n}: \quad \sum_{j=1}^{n} s_j =t  \biggr\}
\end{equation*}
under the covering map $\R^n\to\T^n$,
and has Lebesgue measure zero, so $\Pp_{\xi,i}^t$ is a purely singular measure.
  
Nevertheless, Condition~B guarantees that time averages of transition probabilities have
nontrivial absolutely continuous components.  Specifically, we will establish this for
the resolvent probability kernel $\Qq_{\xi,i}$ defined by 
\begin{equation}
\label{eq:resolvent}
\Qq_{\xi,i}(E\times\{j\}):= \int_{\R_{+}} e^{- t}\,\Pp_{\xi,i}^t(E\times\{j\}) dt.
\end{equation}
The resolvent kernels are useful in the study of invariant distributions due to the following straightforward result.
\begin{lemma}\label{lem:Q-invariance}
If a measure $\mu$ is $(\Pp^t)$-invariant it is also $(\Qq)$-invariant, i.e., 
$\mu=\mu \Qq$, where the convolution $\mu\Qq$ is defined analogously to~\eqref{eq:convolution}. 
\end{lemma}

\begin{theorem}\label{thm:AC-component-2}
If Condition B is satisfied at some point $\xi\in M$, then for any $i\in S$,
the measure $\Qq_{\xi,i}$ defined by~\eqref{eq:resolvent}
has a nonzero absolutely continuous component with respect to Lebesgue measure on $M \times S$. 
\end{theorem}

\smallskip

The convergence of transition probabilities to the invariant measure, provided that it exists, is out of 
the scope of the present paper. In~\cite{Benaim-et-al-Qualitative} it is shown that if Condition~A is
satisfied at a $D$-approachable
point $\xi$, and if $M$ is compact, the transition probabilities converge to a unique invariant measure in total
variation at exponential rate. Our analysis suggests that existence of an invariant measure and Condition~B at a
$D$-approachable point implies only Ces\`aro convergence.

At the heart of our proofs of Theorems~\ref{thm:AC-component-1} and~\ref{thm:AC-component-2} are classical results
from  geometric control theory that can be found in~\cite{Jurdjevic:MR1425878}. The statements we present are derived
from Theorems~$3.1$,~$3.2$, and~$3.3$ in~\cite{Jurdjevic:MR1425878}. Analogous results for the special case of analytic
vector fields on a real analytic manifold are first stated in~\cite[Theorems~$3.1$ and
$3.2$]{Sussmann-Jurdjevic:MR0338882}. In their paper, Sussmann and Jurdjevic were able to build on prior
work~\cite{Chow:MR0001880} by Chow who considered symmetric families of analytic vector fields.
Krener generalized these results to $C^{\infty}$-vector fields in~\cite{Krener:MR0383206}.

Recall that a regular point of a function $f:\mathbb{R}^{m} \to M$ is a point ${\bf t} \in \mathbb{R}^{m}$ such that the differential $Df({\bf t})$ has full rank. 
If $Df({\bf t})$ has deficient rank, ${\bf t}$ is called a critical point of~$f$.

\begin{theorem} \label{thm:Jurdjevic-1} Assume that Condition~A holds at
some $\xi\in M$.
Then:
\begin{enumerate}
 \item For any $i,j\in S$,
there are an integer $m>n$
and  a vector ${\bf i} \in S^{m+1}$ with $i_1=i$ and $i_{m+1}=j$ such that for any $t>0$ the mapping $f_\ibf:\R_{+}^{m} \to M$
defined by
\begin{equation}
\label{eq:flow_for_transition_probability}
f_\ibf(t_1,\ldots,t_{m})=\Phi_{{\bf i}}\biggl(t_1,\ldots,t_{m},t-\sum_{l=1}^{m} t_l, \xi \biggr)
\end{equation}
has a nonempty open set of regular points in the simplex 
\begin{equation*}
\Delta_{t,m}:=\biggl\{(t_1,\ldots,t_m) \in \R_{+}^{m}:\ \sum_{l=1}^{m} t_l < t \biggr\}.
\end{equation*}
\item The interior of $L(\xi)$ is nonempty and dense in $L(\xi)$.
\end{enumerate}
\end{theorem} 
 
\begin{theorem} \label{thm:Jurdjevic-2}
Assume that Condition~B holds at some $\xi\in M$. 
Then:
\begin{enumerate}
 \item For any $i,j\in S$,
there are an integer $m>n$
and  a vector ${\bf i} \in S^{m+1}$ with $i_1=i$ and $i_{m+1}=j$ such that for any $t>0$ the mapping $F_\ibf:\R_{+}^{m+1} \to M$
defined by
\begin{equation*}
F_\ibf(t_1,\ldots,t_{m+1})=\Phi_{{\bf i}}(t_1,\ldots,t_{m+1},\xi)
\end{equation*}
has a nonempty open set of regular points in $\Delta_{t,m+1}$.
\item The interior of $L(\xi)$ is nonempty and dense in $L(\xi)$.
\end{enumerate}
\end{theorem}

Condition~A is stronger than Condition~B, so it is not surprising that the conclusion of
Theorem~\ref{thm:Jurdjevic-1} implies the conclusion of Theorem~\ref{thm:Jurdjevic-2}. We will not prove these
theorems since they are direct consequences of results in~\cite[Chapter 3]{Jurdjevic:MR1425878}.

Theorem~\ref{thm:Jurdjevic-1} shows that under Condition A, we can find a sequence of driving vector fields
such that using that sequence and varying only the switching times we can generate an open set of terminal
positions for any fixed terminal time $t>0$. Moreover, the map assigning the terminal position at time $t$ to the
switching time sequence is regular, i.e., its Jacobian has full rank. We will use this theorem to conclude that, under
this map, the pushforward of
an absolutely continuous measure is also absolutely continuous.

Under Condition B, such regularity for a fixed time $t$ is not guaranteed. However, Theorem~\ref{thm:Jurdjevic-2}
shows that if it is allowed to vary also the terminal time $t$, we still can generate an open set of terminal positions
and the Jacobian of the corresponding map still has full rank. This means that although the
pushforward measures
themselves do not necessarily enjoy
the desired regularity, their averages over terminal times $t$ do, and we will use this argument to study the regularity
of
the resolvent measure of the Markov process under consideration.

The basic idea behind Theorems~\ref{thm:Jurdjevic-1} and~\ref{thm:Jurdjevic-2} is that for a sufficient number of
switches, by perturbing the switching time sequences one can generate perturbations to the
terminal point in all directions.

The first statement of Theorem~\ref{thm:Jurdjevic-2} corresponds to Theorem~$3.1$ in~\cite{Jurdjevic:MR1425878}, which reads as follows: Under the assumptions of Theorem~\ref{thm:Jurdjevic-2}, any neighborhood $U$ of $\xi$ contains points that are normally accessible from $\xi$ at arbitrarily small times. A point $\eta$ in $M$ is called normally accessible from $\xi$ at time $t>0$ if there exist vectors $\ibf \in S^{m+1}$ and $(\hat{t}_1,\ldots,\hat{t}_{m+1}) \in \Delta_{t,m+1}$ such that $F_\ibf(\hat{t}_1,\ldots,\hat{t}_{m+1}) = \eta$ and the differential $DF_\ibf(\hat{t}_1,\ldots,\hat{t}_{m+1})$ has full rank. It's worth pointing out, though, that in~\cite{Jurdjevic:MR1425878} only one sequence~$\ibf$ resulting in $F$ with a regular point is constructed. But since the flow generated by any vector field is a family of diffeomorphisms, and since the set of points satisying Condition~B is open, one can append any indices in front or at the back of that sequence without destroying the desired properties, and thus recover this part of Theorem~\ref{thm:Jurdjevic-2} as we state it.

The fact that the interior of $L(\xi)$ is nonempty and dense in $L(\xi)$ follows from Theorem~$3.2.a$ in~\cite{Jurdjevic:MR1425878}.
Theorem~\ref{thm:Jurdjevic-1} follows from applying Theorem~$3.1$ (\cite{Jurdjevic:MR1425878}) to $\R \times M$ and
vector fields ${\bf 1} \oplus u_i, i\in S$, where
\begin{equation*}
({\bf 1} \oplus u)(r,\xi):= (1,u(\xi)),\quad (r,\xi)\in \R \times M,
\end{equation*}
and ${\bf 1}$ is the unit vector field on $\mathbb{R}$ corresponding to
the derivation $\partial/\partial r$ and identically equal to~$1$ in the natural coordinates on $\R$.

\section{Proof of Theorem~\ref{thm:AC-component-1}} \label{sec:AC-component-1-proof}   

We need to prove that for any $t>0$ and $i \in S$, the measure $\Pp^t_{\xi,i}$ is not singular.

\medskip

For any finite sequence ${\ibf}$ of indices in $S$ with initial index $i$ (we will call these sequences 
admissible), let $C_{{\ibf}}$ be the event that the driving vector fields up to time $t$ appear in the order determined by ${\ibf}$. Since
$\Pp_{\xi,i}(C_{\ibf})>0$ for any admissible $\ibf$
it suffices to find an admissible sequence $\ibf$ such that 
$\Pp_{\xi,i}^t(\cdot| C_{\ibf})$ is not singular.
We claim that this holds true for the the sequence $\ibf$ provided by 
Theorem~\ref{thm:Jurdjevic-1}. According to
Theorem~\ref{thm:Jurdjevic-1}, there is an admissible sequence $\ibf=(i_1,i_2,\ldots,i_{m+1})$ with $i_1 = i$ 
such that the function
$f_\ibf$
has a regular point in $\Delta_{t,m}$. Since the set of regular points of a differentiable function is open in its domain, the function $f_\ibf$ is regular in a nonempty open set $B \subset \Delta_{t,m}$.  

Let $T_1,T_2,\ldots,T_{m+1}$ be independent and exponentially distributed random variables such that $T_j$ has
parameter
$\lambda_{i_j}$ for $1\leq j \leq m+1$.  

On $C_\ibf$ we have $A_t=i_{m+1}$, and the distribution of $X_t$ under $\Pp_{\xi,i}(\cdot | C_{\ibf})$ coincides
with the distribution of $f_\ibf(T_1,\ldots,T_m)$ conditioned on the event
\begin{equation}
\label{eq:event_R}
R=\biggl\{\sum_{j=1}^{m} T_j < t \le \sum_{j=1}^{m+1} T_j \biggr\}.
\end{equation}
The distribution of the random vector $(T_1,\ldots,T_m)$ conditioned on $R$, 
is equivalent to the uniform distribution on the simplex 
\begin{equation*}
\Delta_{t,m} := \biggl\{(t_1,\ldots,t_m) \in \R_{+}^{m}: \quad \sum_{j=1}^{m}t_j < t \biggr\}.
\end{equation*}
Now the theorem directly follows from the following result:

\begin{lemma}\label{lem:pushforward_full_rank} Let $n,m\in\N$, $n\le m$.
Suppose that $B$ and $\Delta$ are nonempty open sets in~$\R^m$, $B\subset\Delta$,
and $M$ is an $n$-dimensional
smooth manifold.
If $f:\Delta\to M$ is 
differentiable on $B$ and all points in $B$ are regular for $f$, then for any
absolutely continuous
probability measure $\mu$ on $\Delta$ satisfying $\mu(B)>0$, its
pushforward
$\mu f^{-1}$ is not singular with respect to~$\lambda^M$.
\end{lemma}

We will prove this lemma only for $M=\R^n$.  Modifying the proof for the general case using
coordinate patches on $M$ amounts only to notational differences. 

We will use the following statement (see, e.g., Proposition~4.4
in~\cite{Davydov-et-al:MR1604537}):
\begin{lemma}\label{lem:pushforward_for_nonzero_det}
 Let  $f:B \to\R^m$ be a Borel function a.e.-differentiable on an open set $B \subset\R^m$
and satisfying  $\lambda^m\{\tbf\in B:\, \det Df(\tbf)=0\}=0$. If $\mu\ll \lambda^m$, then
$\mu f^{-1}\ll\lambda^m$, and
\[
 \frac{d(\mu f^{-1})}{d\lambda^m}(\sbf)=\sum_{\tbf\in B: f(\tbf)=\sbf} |\det Df(\tbf)|^{-1}\frac{d\mu}{d\lambda^m}(\tbf)
\]

\end{lemma}

\bpf[Proof of Lemma~\ref{lem:pushforward_full_rank}]
We can find an open  set $B'\subset B$ such that $\mu(B')>0$ and there are $n$
columns of $Df(\tbf)$ (without loss of generality, first~$n$ columns) such that for any $\tbf \in B'$ they are linearly
independent.
For $\rho:B' \to \R^n \times \mathbb{R}^{m-n}$ defined by
\[
\rho:\tbf=(t_1,\ldots,t_m) \mapsto (f(\tbf), t_{n+1},\ldots,t_m),
\]
and any $\tbf\in B'$, we have
$
\det D\rho(\tbf ) \neq 0. 
$
Therefore, by Lemma~\ref{lem:pushforward_for_nonzero_det}, the pushforward of
the restriction of $\mu$ to $B'$ under $\rho$
is a positive absolutely continuous measure on $\R^n \times \mathbb{R}^{m-n}$. Integrating over $\mathbb{R}^{m-n}$, we
obtain that the pushforward of the restriction of $\mu$ to $B'$ under $f$
is a positive absolutely continuous measure on $\R^n$, and the proof is complete.
\epf

\section{Proof of Theorem~\ref{thm:AC-component-2}}  \label{sec:AC-component-2-proof}

 We need to show that $\Qq_{\xi,i}$ is not a singular measure. The proof is based on
Theorem~\ref{thm:Jurdjevic-2}.

For the $S$-valued process $A$ we denote by $I_t(A)$ the sequence of states visited by $A$ between $0$
and $t$.
For any $m\in\N$ and any sequence $\ibf\in S^m$, we can introduce an auxiliary measure $\Qq_{\xi,i,\ibf}$ on $M$ by 
\[
\Qq_{\xi,i,\ibf}(B)=\int_{R_{+}} e^{-t}\,\Pp_{\xi,i}\{X_t\in B\ \text{and}\ I_t(A)=\ibf\}\,  dt,\quad B\in\Bc(M).
\]
Since 
\begin{equation}
\label{eq:Qq-via-I_t}
\Qq_{\xi,i}(B \times \{j\})=\sum_{m}\sum_{\ibf=(i,i_2,\ldots,i_{m-1},j)\in S^m}\Qq_{\xi,i,\ibf}(B), 
\end{equation}
it is sufficient to find $\ibf=(i_1,\ldots,i_m)$ with $i_1=i$ such that $\Qq_{\xi,i,\ibf}(M)>0$ and
\[
\wb \Qq_{\xi,i,\ibf}(\cdot)= \frac{\Qq_{\xi,i,\ibf}(\cdot)}{\Qq_{\xi,i,\ibf}(M)} 
\]
 is a nonsingular probability measure. To apply Lemma~\ref{lem:pushforward_full_rank}, we need to represent
$\wb \Qq_{\xi,i,\ibf}$ as the pushforward of a measure, equivalent to Lebesgue measure, under a smooth map with a 
nonempty set of regular points.

Since Condition~B holds at $\xi$, Theorem~\ref{thm:Jurdjevic-2} yields an integer $m>n$
and a sequence $\ibf=(i_1,i_2,\ldots,i_{m+1})$ with $i_1=i$, such that the function $F_\ibf:\R_{+}^{m+1} \to
M$ defined by
\begin{equation*}
F_\ibf(\tbf)= \Phi_{\ibf}(\tbf,\xi)
\end{equation*}
has a regular point.
For this $\ibf$ provided by Theorem~\ref{thm:Jurdjevic-2}, $\wb \Qq_{\xi,i,\ibf}$ is 
the distribution of $\Phi_{{\bf i}} \left(T_1,\ldots,T_m,T-\sum_{j=1}^{m} T_j, \xi \right)$ conditioned on
\begin{equation}
R= \biggl\{\sum_{j=1}^{m} T_j <T \leq \sum_{j=1}^{m+1} T_j \biggr\},
\label{eq:event_R_2}
\end{equation}
where $T_1,\ldots,T_{m+1}$, and $T$ are independent random variables 
exponentially distributed with parameters
$\lambda_{i_1},\ldots,\lambda_{i_{m+1}}$, and $1$, respectively.

Since the joint distribution of $T_1,\ldots,T_{m+1},T$ is equivalent to Lebesgue measure and since event $R$ has positive probability,
the distribution $\mu$ of $T_1,\ldots,T_m,T$ conditioned on $R$
%
induces a measure on
\[
\Delta= \biggl\{(t_1,\ldots,t_m,t) \in \R_{+}^{m+1}: \quad \sum_{j=1}^{m} t_j < t \biggr\}
\]
that is equivalent to Lebesgue measure. The regularity of $F_\ibf$ guaranteed by Theorem~\ref{thm:Jurdjevic-2} implies
that
the function $f_\ibf:\Delta\to M$ defined 
by
\begin{equation*}
 f_\ibf(t_1,\ldots,t_m,t)= F_\ibf\biggl(t_1,\ldots,t_m,t-\sum_{j=1}^{m} t_j \biggr)
\end{equation*}
has a nonempty open set of regular points in~$\Delta$, and the proof is completed by an application of
Lemma~\ref{lem:pushforward_full_rank}, since $\wb \Qq_{\xi,i,\ibf}$ is the pushforward of 
$\mu$ under $f_\ibf$.

\section{Absolute continuity of ergodic invariant measures}\label{sec:AC-of-ergodic-measures}

According to the Ergodic Decomposition Theorem, all invariant measures for a Markov semigroup can be
represented as
convex combinations of ergodic ones (see, e.g., \cite[Theorem~1.7]{Hairer:ergodicity-lectures}). We will use this to derive Theorem~\ref{thm:uniqueness} from absolute continuity of ergodic invariant distributions.

To define ergodicity, we need to recall the notion of $\mu$-invariant sets. Let $\mu$ be an invariant measure for the
Markov semigroup $(\Pp^t)$. We say that a set $A\in\Bc(M)
\otimes \Pc(S)$ is $\mu$-invariant if for every $t\ge 0$, $\Pp^t_{\xi,i}(A)=1$ for $\mu$-almost every
$(\xi,i)\in A$. An invariant measure $\mu$ is called ergodic if for every $\mu$-invariant set $A$, either  $\mu(A)=1$ or
$\mu(A)=0$.

The following is a basic result on systems with Markov switchings that does not use Conditions A or B.

\begin{theorem}\label{thm:AC-or-S}
If $\mu$ is $(\Pp^t)$-invariant and ergodic then it is either absolutely continuous or singular. 
\end{theorem}
\bpf
Consider the Lebesgue decomposition $\mu=\mu_{ac}+\mu_{s}$, where $\mu_{ac}$ is absolutely
continuous and $\mu_{s}$ is singular with respect to Lebesgue measure.  Let us show that both
$\mu_{ac}$ and $\mu_{s}$ are invariant.

For any $t>0$, using the
invariance of $\mu$, we can write
\begin{align}
\mu_{ac}+\mu_{s} &= \mu=\mu \Pp^t
\label{eq:decomposition_in_ac_and_s}
=
\mu_{ac} \Pp^t+\mu_s\Pp^t
=\sum_{j=1}^{k} \nu_j  + \mu_s\Pp^t,
\end{align}
where 
\begin{equation}
\label{eq:nu_j}
\nu_{j}(\cdot)= \int_M \Pp^t_{\xi,j} (\cdot)\mu_{ac}(d\xi\times\{j\}),\quad j\in S.
\end{equation}

We claim that the measures $\nu_j$, $j\in S$,
are absolutely continuous. To see this we
check that for any sequence $\ibf=(i_1,\ldots,i_{m+1})$ with $i_1=j$, the measure $\nu_\ibf$ defined by
\begin{align}\notag
 \nu_{\ibf}(E)&= \int_M \Pp_{\xi,j} (X_t\in E\, |C_\ibf)\mu_{ac}(d\xi\times\{j\})\\
&=\int_M \Pp\biggl(\Phi_{{\bf i}}\biggl(T_1,\ldots,T_m,t-\sum_{l=1}^{m} T_l,\xi \biggr)\in
E\, \biggr|\, R\biggr)\mu_{ac}(d\xi\times\{j\})\label{eq:one_of_continuous_components_disintegrated}
\end{align}
 is absolutely continuous (here we use the notation introduced in Section~\ref{sec:AC-component-1-proof}). Suppose
$\lambda^M(E)=0$. For fixed $T_1,\ldots,T_m,T_{m+1}$,  the map $\Phi_\ibf$ is a diffeomorphism in $\xi$. Therefore, on
event $R$ introduced in~\eqref{eq:event_R}, we have
\[
 \mu_{ac}\biggl(\xi \times \{j\}:\ \Phi_{{\bf i}}\biggl(T_1,\ldots,T_m,t-\sum_{l=1}^{m} T_l,\xi \biggr)\in
E  \biggr)=0,
\]
and $\nu_\ibf(E)=0$ follows from disintegrating the right side
of~\eqref{eq:one_of_continuous_components_disintegrated} and changing the order of integration.

Now, using~\eqref{eq:decomposition_in_ac_and_s} and the absolute continuity of $\nu_j$, $j\in S$, we can write
\begin{equation}
\label{eq:decomposition_of_ac_part}
\mu_{ac}= 
\sum_{j=1}^{k} \nu_j  + (\mu_s\Pp^t)_{ac}\ .
\end{equation}
Since  $\Pp^t_{\xi,j}(M\times S)=1$ for all $\xi$ and $j$, \eqref{eq:nu_j} implies
$\sum_{j=1}^{k} \nu_j(M\times S)= \mu_{ac}(M\times S)$. Therefore, applying~\eqref{eq:decomposition_of_ac_part}
to $M \times S$, we obtain that the absolutely continuous
component of the measure $\mu_s\Pp^t$ is zero. In other
words, $\mu_s\Pp^t$ is singular, and from~\eqref{eq:decomposition_in_ac_and_s} and the absolute
continuity of $\nu_j$, $j\in S$, we obtain
\begin{equation}
\label{eq:invariance_of_singular_component}
\mu_s = \mu_s\Pp^t.
\end{equation}                        
In other words, $\mu_s$ is invariant for $(\Pp^t)$. It follows from~\eqref{eq:decomposition_in_ac_and_s}
that $\mu_{ac}$ is also invariant. Since $\mu$ is ergodic, it cannot be represented as a sum of two nontrivial
invariant measures. This means that either $\mu=\mu_{ac}$ or $\mu=\mu_s$.
\epf

We endow the state space $S$ with the discrete topology and recall that a point $(\xi,i)\in M\times S$ is
contained in the support of a measure if and
only if the measure of every open neighborhood of $(\xi,i)$ is positive.  
 
\begin{theorem}\label{thm:absolute-continuity-for-ergodic}
Let $\mu$ be an ergodic invariant measure for $(\Pp^t)$. Assume that the support of $\mu$ contains a point
$(\eta,i)$ such that Condition~B holds at~$\eta$. 
Then, $\mu$ is absolutely continuous
with
respect to Lebesgue measure on $M\times S$.  
\end{theorem}

We will need several auxiliary statements.

\begin{lemma} \label{lem:support}
Let $\nu$ be a finite Borel measure on $M \times S$ with support $K$. If $U$ is any open set in $M \times
S$ whose intersection with $K$ is nonempty, we have 
\begin{equation*}
\nu(U \cap K)>0.
\end{equation*}
\end{lemma} 
 
\bpf
Assume that $\nu(U \cap K)=0$. The complement of the support $K$ has measure
zero. Therefore
\begin{equation*}
\nu(U) = \nu(U \cap K) + \nu(U \cap K^{c})=0.
\end{equation*}
Thus, $U^{c}$ is a closed subset of $M \times S$ whose complement has measure zero. From the
definition of the support, we obtain $K \subset U^{c}$. But then, $U \cap K$ must be empty, a contradiction. 
\epf

\bpf[Proof of Theorem~\ref{thm:absolute-continuity-for-ergodic}]  According to Theorem~\ref{thm:AC-or-S} we need to show
that $\mu$ is not singular. If $\mu$ is singular,
it is entirely supported on a zero Lebesgue measure set $G\subset M\times S$, so
$\mu(G^c)=0$. Since $\mu$ is $(\Pp^t)$-invariant, it is also $\Qq$-invariant. Therefore,
$\mu(G^c)=\mu \Qq(G^c)$, and we see that $\mu(V)=0$ where
\[
 V=\{(\xi,j)\in M \times S:\ \Qq_{\xi,j} (G^{c})>0\}.
\]
Let $U$ be the set of points $\xi \in M$ where Condition~B holds.
Due to Theorem~\ref{thm:AC-component-2}, $U\times S \subset V$, and
we conclude that $\mu(U\times S)=0$.

Recall that $U$ is an open subset of $M$, and $(U\times S)\cap \supp \mu\ne\emptyset$ by assumption.
Lemma~\ref{lem:support} implies that  $\mu((U\times S)\cap \supp \mu)>0$.  The contradiction with $\mu(U\times
S)=0$ completes the proof.
\epf

Of course, if one replaces Condition~B in Theorem~\ref{thm:absolute-continuity-for-ergodic} with 
the stronger Condition~A the resulting statement holds automatically, but one can give a proof that does not involve
the resolvent $Q$:

\begin{theorem}\label{thm:absolute-continuity-for-ergodic-under-condition-A}
Let $\mu$ be an ergodic invariant measure for $(\Pp^t)$. Assume that the support of $\mu$ contains a point
$(\eta,i)$ such that Condition~A holds at~$\eta$.
Then, $\mu$ is absolutely continuous
with
respect to the Lebesgue measure on $M\times S$.  
\end{theorem}

\bpf  According to Theorem~\ref{thm:AC-or-S} we need to show that $\mu$ is not singular. If $\mu$ is singular,
it is entirely supported on a zero Lebesgue measure set $G\subset M\times S$, so
$\mu(G^c)=0$. Since $\mu(G^c)=\mu \Pp^t(G^c)$, we see that $\mu(V)=0$ where
\[
 V=\{(\xi,j)\in M \times S:\ \Pp^t_{\xi,j} (G^{c})>0\}.
\]
Let $U$ be the set of points $\xi \in M$ where Condition~A holds.
Due to Theorem~\ref{thm:AC-component-1}, $U\times S \subset V$, and
 we conclude that $\mu(U\times S)=0$.

Recall that $U$ is an open subset of $M$, and $(U\times S)\cap \supp \mu\ne\emptyset$ by assumption.
Lemma~\ref{lem:support} implies that  $\mu((U\times S)\cap \supp \mu)>0$.  The contradiction with $\mu(U\times
S)=0$ completes the proof.
\epf

\section{Proof of Theorem~\ref{thm:uniqueness}} \label{sec:uniqueness-proof}

First, we establish two properties of the set $E= L\cap U$, where
$U$ is the open set of points satisfying Condition B.
\begin{lemma} \label{lem:non-empty-interior}
The set $E$ has nonempty interior.
\end{lemma}

\bpf  By assumption,  $\xi\in E$ , so $U\ne\emptyset$ and $L(\xi) \cap U\ne \emptyset$ by continuity of the vector fields in $D$. Since $\xi\in U$, Theorem~\ref{thm:Jurdjevic-2} implies that $L(\xi)$
has nonempty interior that is dense in $L(\xi)$. Therefore, the set
\begin{equation*}
V=L(\xi)^{\circ} \cap U
\end{equation*}
is nonempty and open. Clearly, $V\subset U$, and it remains to prove that $L(\xi)^{\circ}\subset L$.
In fact, we even have that $L(\xi) \subset L$. To see that, let us fix any  $\zeta\in L(\xi)$, $\eta\in M$, 
and prove
that $\zeta\in \overline{L(\eta)}$.
Since $\zeta\in L(\xi)$, we have  
\begin{equation*}
\zeta = \Phi_{{\bf i}}({\bf t},\xi)
\end{equation*}
for some index sequence $\ibf$ and some time sequence $\tbf$. Let us fix a neighborhood $W$ of $\zeta$. 
Since the mapping $x \mapsto \Phi_{{\bf i}}({\bf t},x)$ is continuous, the inverse image of $W$ under this map is an
open neighborhood of $\xi$. Since $\xi$ is $D$-approachable from $\eta$, this
open neighborhood of $\xi$ contains a point $D$-reachable from $\eta$. Hence, $W$ contains a point that is $D$-reachable
from $\eta$.
\epf

As an immediate corollary of Lemma~\ref{lem:non-empty-interior}, the set $L$ has nonempty interior.
\begin{lemma} \label{lem:positive-measure}
Suppose $\mu$ is an invariant measure for  $(\Pp^t)$. If $G$ is a nonempty open subset of~$L$ and $j\in
S$, then $\mu(G\times\{j\})>0$.
\end{lemma}      
  
\bpf Let us assume that $\mu(G\times\{j\})=0$. Since
$\mu$ is $(\Pp^t)$-invariant, it is also $\Qq$-invariant, and we have
\begin{equation*}
0 = \mu(G\times\{j\}) = \sum_{i=1}^{k} \int_M \Qq_{\eta,i}(G\times\{j\}) \mu(d\eta\times\{i\}).
\end{equation*}
For all $i\in S$ and $\mu(\cdot \times\{i\})$-almost every $\eta\in M$, we thus obtain
\begin{equation}
\label{eq:Q-0}
\Qq_{\eta,i}(G\times\{j\}) = 0.
\end{equation}
Let us choose $\eta$ such that~\eqref{eq:Q-0} holds true.

By assumption, we have $G\subset L\subset\overline{L(\eta)}$. Since $G$ is open,  $G \cap L(\eta)\ne
\emptyset$. So there exist  a sequence $\ibf=(i,i_2,\ldots,i_m,j)$ and an interswitching time vector
$\tbf=(t_1,\ldots,t_m,t_{m+1})$ such that  $\Phi_{{\bf i}}({\bf t}, \eta)\in G$. By continuity of $\Phi_{\bf i}$ there
is a neighborhood $W$ of $\tbf$ in $\R^{m+1}_+$
such that $\Phi_{{\bf i}}(\sbf, \eta)\in G$ for all $\sbf\in W$.  Denoting $s=s_1+\ldots+s_{m+1}$ and using the
representation of $\Pp^s_{\eta,i}(\cdot|C_{\ibf})$ via exponentially distributed times
that we used in the proof of Theorem~\ref{thm:AC-component-1}, we conclude that   $\Pp^s_{\eta,i}(G \times \{j\})>0$ for 
$s$ sufficiently close to $t=t_1+\ldots+t_{m+1}$. Therefore, $\Qq_{\eta,i}(G \times \{j\})>0$ contradicting \eqref{eq:Q-0}.
\epf

\bigskip

\bpf[Proof of Theorem~\ref{thm:uniqueness}]
As a consequence of Birkhoff's ergodic theorem, any invariant measure can be written as a convex combination
of ergodic invariant measures, see, e.g.,  Theorem~1.7 in~\cite{Hairer:ergodicity-lectures}. Therefore, it
suffices to
show absolute continuity and uniqueness of an ergodic invariant measure. 

Let us begin by deriving absolute continuity. If $\mu$ is an
ergodic invariant measure that satisfies the assumptions of Theorem~\ref{thm:uniqueness} then, 
due to Theorem~\ref{thm:absolute-continuity-for-ergodic}, it suffices to show that
$L \subset \supp \mu$. Let $\xi \in L$, and let $U$ be a neighborhood of $\xi$ in $M$, $j \in S$. By
Lemma~\ref{lem:positive-measure}, we have $\mu(U \times \{j\})>0$, hence $\xi \in \supp \mu$.   

In order to prove uniqueness of the ergodic invariant measure, let us assume that $\mu_1$ and $\mu_2$ are two distinct ergodic invariant probability measures. Birkhoff's ergodic
theorem then implies that $\mu_1$ and $\mu_2$ are mutually singular. Hence, the set $M \times S$ can be partitioned into
two disjoint subsets $H_1$ and $H_2$ with
$\mu_1(H_2) = \mu_2(H_1) = 0.$
The two sets can be represented as 
\[H_\alpha=\bigcup_{j=1}^{k} M_{\alpha,j} \times \{j\},\quad \alpha=1,2,\]
for some measurable sets $M_{\alpha,j}$, $j\in S$, $\alpha=1,2$. For all $\alpha$ and $j$,
\begin{equation*}
\mu_\alpha( M_{\alpha,j} \times \{j\})=\mu_\alpha( M \times \{j\})>0,
\end{equation*}
since the left side is a stationary distribution for the Markov chain on $S$ and by our assumptions, transitions
between all states happen with positive probability.

Fix a $j$ in $S$. By virtue of Lemma~\ref{lem:non-empty-interior}, the set $E^{\circ}$ is nonempty.
According to Lemma~\ref{lem:positive-measure}, for all $j \in S$ we have $\mu_1(E^\circ\times\{j\})>0$.
 Since $\mu_1(M_{2,j}\times\{j\}) =0$, we deduce that $\mu_1(E_1\times\{j\}) > 0$, where $E_1=E^\circ \cap
M_{1,j}$.

The measure $\mu_1$ is $(\Pp^t)$-invariant, hence it is also $\Qq$-invariant, and we have
\begin{equation}
\label{eq:invariance-for-M_2_j}
0=\mu_1(M_{2,j}\times\{j\}) \geq \int_{E_1} \Qq_{\eta,j}(M_{2,j}\times \{j\}) \mu_1(d\eta\times\{j\}).
\end{equation}
Since $\mu_1(E_1\times\{j\})>0$, it suffices to show that $\Qq_{\eta,j}(M_{2,j}\times\{j\})>0$  for all $\eta \in
E_1$, to obtain a contradiction with~\eqref{eq:invariance-for-M_2_j}.

Since $\eta$ satisfies Condition~B,
Theorem~\ref{thm:Jurdjevic-2} guarantees that there exist an integer $m>n$ and a vector ${\bf i} =
(j,i_2,\ldots,i_m,j)$ such
that the function $f:\mathbb{R}^{m+1}_+ \to M$ defined by
\begin{equation*}
 f(\bf t) =\Phi_{\bf i}({\bf t},\eta)
\label{eq:F_ibf}
\end{equation*}
has an open set $O$ of regular points such that for all $t>0$,
\begin{equation*}
\{\tbf=(t_1,\ldots,t_{m+1})\in O:\ t_1+\ldots+t_{m+1}<t\}\ne \emptyset.
\end{equation*}
Therefore, the map $F$ defined by
\[
 F(t_1,\ldots,t_{m+1},t)= f\biggl(t_1,\ldots,t_m,t-\sum_{l=1}^{m}t_l\biggr)
\]
on
\[
 \Delta=\biggl\{(t_1,\ldots,t_{m+1},t)\in\R_{+}^{m+2}:\ \sum_{l=1}^{m}t_l<t<\sum_{l=1}^{m+1}t_l\biggr\}
\]
has an open set $V\subset\Delta$ of regular points such that 
\begin{equation}
\{\tbf=(t_1,\ldots,t_{m+1},t)\in V:t<s\}\ne \emptyset,\quad s>0.
\label{eq:0_is_limiting_point}
\end{equation}

Using the representation of $\Qq$ via~\eqref{eq:Qq-via-I_t} and the family of exponentially distributed times
$T_1,\ldots,T_{m+1},T$, we obtain that it is sufficient to prove that
\begin{equation}
\Pp\left\{F(T_1,\ldots,T_{m+1},T) \in M_{2,j}|\ R \right\} > 0,
\label{eq:sufficient_for_contradiction}
\end{equation}
where $R$ was introduced in~\eqref{eq:event_R_2}.

 Since $E^\circ$ is an open set containing
$\eta$, and $F(V)$ is an open set such that $\eta\in \overline{F(V)}$ (due 
to~\eqref{eq:0_is_limiting_point} and continuity of $F$ at $0$), we obtain that
$G= E^{\circ} \cap F(V)$ is also a nonempty open set.

Let us choose a vector ${\bf r}\in V$ such that $F({\bf r}) \in E^{\circ}$. Since $\rbf$ is a regular point for
$F$, we see that 
for an arbitrary choice of local smooth coordinates around $\rbf$,
there are $n$ independent columns of the matrix $DF(\sbf)$ for $\sbf$
in a small neighborhood of $\rbf$. Without loss of generality we
can assume that these are the first $n$ columns. 
Then the map $\rho:\R^{m+2}\to M\times\R^{m+2-n}$ defined by
\[
 \rho(s_1,\ldots,s_{m+1},s)=(F(s_1,\ldots,s_{m+1},s),s_{n+1},\ldots ,s_{m+1},s)
\]
has nonzero Jacobian in that neighborhood. So we can choose an open set $W_V$ containing~$\rbf$ so that $\rho$ is a
diffeomorphism between $W_V$ and 
$W_G\times W_{n-m-2}$, where $W_G\subset G$ and $W_{m+2-n}\subset \R_{+}^{m+2-n}$ are some open sets.

The set $W_G$ is an open subset of $L$. It is also not empty since it contains~$F({\bf r})$.
Lemma~\ref{lem:positive-measure} implies that $\mu_2( W_G\times\{j\})>0$.  Since
$\mu_2(M_{2,j}^c\times\{j\})=0$, we conclude that
$\mu_2(J\times\{j\})>0$ where $J=M_{2,j} \cap W_G.$ Since $\mu_2$ is an ergodic measure, it is absolutely continuous, so
\begin{equation}
\label{eq:J-Lebesgue-positive}
\lambda^M(J)>0. 
\end{equation}

Since $J\subset M_{2,j}$, the desired inequality \eqref{eq:sufficient_for_contradiction} will follow from
\begin{equation}
\Pp\left\{F(T_1,\ldots,T_{m+1},T) \in J|\ R \right\} > 0.
\label{eq:sufficient_for_contradiction-2}
\end{equation}

Since the joint distribution of $T_1,\ldots,T_m,T_{m+1},T$ is equivalent to the Lebesgue measure on $\Delta$, 
Lemma~\ref{lem:pushforward_for_nonzero_det} implies that $\rho(T_1,\ldots,T_{m+1},T)$ has positive density
almost everywhere in $W_G\times W_{m+2-n}$. Integrating over $W_{m+2-n}$, we see that \\ $F(T_1,\ldots,T_{m+1},T)$ has
positive density almost everywhere in $W_G$. Now~\eqref{eq:sufficient_for_contradiction-2} follows from 
\eqref{eq:J-Lebesgue-positive}.\epf

\bigskip

Of course, Theorem~\ref{thm:uniqueness} remains true if one replaces Condition~B by the stronger Condition~A.
However, under that condition one can prove this result without referring to the resolvent $\Qq$. Namely, one can
use the regularity of transition probabilities established in Theorem~\ref{thm:AC-component-1} (which is stronger than
the regularity established in Theorem~\ref{thm:AC-component-2}),
and invoke Theorems~\ref{thm:Jurdjevic-1} and~\ref{thm:absolute-continuity-for-ergodic-under-condition-A} instead of
Theorems~\ref{thm:Jurdjevic-2} and~\ref{thm:absolute-continuity-for-ergodic}.

\section{Examples} \label{sec:examples}

In this section, we apply Theorem~\ref{thm:uniqueness} to two concrete switching systems. In the first example, we have a closer look at the system on the $n$-dimensional torus $\T^n=\R^n/\Z^n$ that was introduced in Section~\ref{sec:absolute-continuity-invariant}. In the second example, we switch between two Lorenz vector fields with different parameter sets. For both systems, uniqueness of the invariant measure is derived from Theorem~\ref{thm:uniqueness}. For the system on $\T^n$, we point out the invariant measure explicitly.

\bigskip

Let $M$ be the $n$-dimensional torus $\T^n$, and let $D=\{u_1,\ldots,u_n\}$ be the standard basis of $\mathbb{R}^{n}$.
 We assume for simplicity that the parameter $\lambda$ of the exponential time between any two switches is independent
of the current state, and that we have a uniform probability of switching between any two states. In
Section~\ref{sec:absolute-continuity-invariant} we implicitly argued that Condition~A does not hold for this system: If
Condition~A was satisfied at some point $\xi \in \T^n$, the transition probability measures $\Pp_{\xi,i}^t$ would not be
singular with respect to Lebesgue measure, according to Theorem~\ref{thm:AC-component-1}. However, as pointed out in
Section~\ref{sec:absolute-continuity-invariant}, the measures $\Pp_{\xi,i}^t$ are purely singular. 

It is also instructive to show directly why Condition~A does not hold. As all the vector fields in $D$ are constant, the derived algebra $\Ic'(D)$ contains only the zero vector field. Thus, for any $\xi \in \T^n$,
\begin{equation*}
 \Ic_{0}(D)(\xi) = \biggl\{\sum_{i=1}^{n} \lambda_{i} u_{i}:\ \sum_{i=1}^{n} \lambda_{i} = 0\biggr\}. 
\end{equation*}
Due to the constraint $\sum_{i=1}^{n} \lambda_{i} = 0$, the algebra $\Ic_{0}(D)(\xi)$ does not have full
dimension, so Condition~A is violated at every point in $\T^n$.

On the other hand, Condition~B is clearly satisfied at any point $\xi \in \T^n$, as the standard basis of $\R^n$ applied to $\xi$ yields a full-dimensional set of vectors in the tangent space. Also note that any point in $\T^n$ is $D$-reachable from any other point. Therefore, Theorem~\ref{thm:uniqueness} guarantees that the associated Markov semigroup has a unique invariant measure, provided that such a measure exists. In this elementary example, it is possible to point out the invariant measure explicitly. For Borel sets $E \subset \T^n$ and states $i \in S$, it is given by
\begin{equation*}
 \mu(E \times \{i\}) = \frac{1}{n} \cdot \lambda(E).
\end{equation*}
Here, $\lambda$ denotes Lebesgue measure on $\T^n$. 

\bigskip

The second example provides a situation where (i) the number of vector fields
in~$D$ is less than the
dimension of the manifold $M$, and (ii) each individual vector field in~$D$ gives rise to dynamics with a
strange attractor and no absolutely continuous invariant measures, but (iii) the switched system has a unique invariant
measure and it is absolutely continuous.

Namely, we consider switching between two Lorenz vector fields with different parameter values. A Lorenz
vector field is a vector field defined in $\R^3$, of the form
\begin{equation*}
 u(x,y,z) = \begin{pmatrix}
             \sigma \cdot (y-x) \\
             rx-y-xz \\
             xy-bz
            \end{pmatrix},
\end{equation*}
 where $\sigma$, $r$ and $b$ are physical parameters. Let the set $D$ contain exactly two Lorenz vector fields $u_1$ and
$u_2$ such that $u_1$ has Rayleigh number $r=r_1=28$ and $u_2$ has a Rayleigh number $r=r_2$ different from, but close
to, $28$. We
assume for both vector fields that $\sigma=10$ and that $b=\tfrac{8}{3}$, which is the classical parameter choice for
the Lorenz system. In~\cite{Tucker:MR1701385}, Tucker shows that the Lorenz system with parameters $\sigma=10$, $r=r_1$
and $b=\tfrac{8}{3}$, corresponding to vector field $u_1$, admits a robust strange attractor $\Lambda$ as well as a
unique SRB-measure supported on $\Lambda$. Robustness implies that the dynamical structure of the system remains intact
under small parameter changes, so the dynamics induced by $u_2$ share these features if $r_2$ is sufficiently close to
$r_1$. Moreover, the SRB-measure on $\Lambda$ satisfies a dissipative ergodic theorem, see
e.g.~\cite[Section~5.1]{Bunimovich-Sinai:MR755521}. It follows that any point $\xi\in\Lambda$ is
$\{u_1\}$-approachable (and thus $D$-approachable) from every point in a set $S_\xi\subset\R^3$ with zero Lebesgue
measure complement.

Assisted by a computer algebra system, we checked that Condition~A is satisfied for this system at any point in $\R^3$
that does not lie on the $z$-axis. Since the $z$-axis is invariant under the flows of both vector fields, we disregard
it and set $M$ to be $\R^3$ without points on the $z$-axis. With this provision, every point on the attractor $\Lambda$
is $D$-approachable from any point in $M$: 

Consider a point $\xi \in \Lambda$ and a point $\eta \in M$. By Theorem~\ref{thm:Jurdjevic-1}, there is a nonempty open
set of $D$-reachable points from $\eta$ (recall that Condition~A holds at any point in $M$). And since this open set has
positive Lebesgue measure, it contains a point belonging to~$S_\xi$.
Hence, $\xi$ is $D$-approachable from $\eta$.  As in the first example, uniqueness and absolute continuity of an invariant measure follow now from Theorem~\ref{thm:uniqueness}.

\begin{remark}   \rm 
In~\cite{Bakhtin-Hurth}, the version of this article that was published in Nonlinearity, we erroneously claimed that existence of an invariant measure for the Lorenz switching system described above followed from the fact that one can construct a common Lyapunov function for the vector fields $u_1$ and $u_2$. Edouard Strickler pointed out to us that existence of such a Lyapunov function only implies existence of an invariant measure for the switching system considered on $\R^3 \times \{1,2\}$, without removing the $z$-axis. Such an invariant measure, however, trivially exists: One can simply take $\delta_0 \otimes \nu$, where $\delta_0$ is the Dirac measure at the origin and where $\nu$ is the unique stationary distribution for the Markov chain $(A_t)_{t \geq 0}$.  In~\cite[Proposition~3.1]{Strickler}, Strickler shows as a corollary to his theory on invariant measures for switching systems with a common equilibrium on a shared invariant face that even the Lorenz switching system on $M \times \{1,2\}$ admits an invariant measure. We would like to thank him for closing this gap. 
\end{remark}

\bibliographystyle{alpha}
\bibliography{switch}

\end{document}